\documentclass[11pt]{amsart}

\usepackage{amsmath, amsthm, amsfonts}
\newtheorem{proposition}{Proposition}
\newtheorem{theorem}{Theorem}

\newtheorem{lemma}{Lemma}

\newcommand{\ze}{\mathbb{Z}}

\author[Milton Jara, Gregorio Moreno and Alejandro F. Ram\'{\i}rez]{Milton Jara, Gregorio Moreno$^1$ and Alejandro F. Ram\'{\i}rez$^1$}

\address[Milton Jara]{\noindent IMPA, Estrada
Dona Castorina 110, CEP 22460 Rio de Janeiro, RJ, Brazil
\newline
e-mail:  \rm \texttt{monets@impa.br}}

\address[Gregorio Moreno]{Facultad de Matem\'aticas\\
Pontificia Universidad Cat\'olica de Chile\\
Vicu\~na Mackenna 4860, Macul\\
Santiago 6904441, Chile
\newline
e-mail:  \rm \texttt{gmoreno@puc.cl}}

\address[Alejandro F. Ram\'\i rez]{Facultad de Matem\'aticas\\
Pontificia Universidad Cat\'olica de Chile\\
Vicu\~na Mackenna 4860, Macul\\
Santiago 6904441, Chile
\newline
e-mail:  \rm \texttt{aramirez@mat.puc.cl}}

\thanks{ AMS 2000 {\it subject classifications}. Primary  60F17, 82C22, 82C41;
 secondary 82B24, 60K35, 60G99.}


\thanks{$^1$Partially supported by Fondo Nacional de Desarrollo Cient\'\i fico
y Tecnol\'ogico grant 1060738}

\date{}

\title[Front propagation in an exclusion dynamics]{Front propagation
in an exclusion one-dimensional reactive dynamics}

\begin{document}

\begin{abstract} We consider an exclusion process
representing a  reactive dynamics of a
pulled front on the integer lattice, describing the dynamics of first class $X$ particles moving as a
simple symmetric exclusion process, and static second class $Y$ particles.
When an $X$ particle  jumps to a site
with a $Y$ particle,  their position is
intechanged and the $Y$ particle becomes an $X$
one. Initially, there is an arbitrary configuration of
$X$ particles at sites $\ldots, -1,0$,  and  $Y$ particles only at sites
$1,2,\ldots$, with a product Bernoulli law of parameter $\rho,
0<\rho<1$. We prove
a law of large numbers and a  central limit theorem
for the front defined by the right-most visited site of the  $X$ particles at time $t$.
 These results corroborate  Monte-Carlo
simulations performed  in
a similar context.
 We also prove that the law of the $X$ particles
 as seen from the front
converges  to a unique invariant measure.
The proofs use regeneration times: we
present a direct way
 to define them within this context.

\end{abstract}

\maketitle

\section{Introduction}

Few mathematical results
exist about  microscopic models of pulled front
propagation 
 representing non-equilibrium pattern formation in chemical reactions,
or physical or biological phenomena
(see \cite{p} and \cite{s} for a review of the physical literature).
One-dimensional interacting particle systems
which are microscopic versions of the Fisher-Kolmogorov-Petrovsky-Piscunov
equation,
representing  evolutionary phenomena in genetics,
have been studied in \cite{bcdls} and \cite{dp}.
More recently, in several works (see \cite{ks},  \cite{m} and \cite{cqr2}),  systems of branching interacting
random walks on the  lattice were studied, which 
consider the spatial ordering of particles: two types of
particles perform independent continuous time simple random walk movement: 
 $X$
particles which  jump at rate $D_X$ and  $Y$ particles at
rate $D_Y$. Upon contact with an $X$ particle, a $Y$ particle
becomes $X$. In \cite{cqr1} the one-dimensional case where $D_Y=0$,
representing the combustion of a propellant towards a stationary state,
was analyzed.
There, the particles are symmetric random walks
 with initially  one $Y$ particle
at each site $1,2,\ldots$ and an arbitrary configuration
of $X$ particles at $-2,-1,0$ with a finite $l_1$ norm
with a certain exponential weight. If we call $r_t$ the position
of the right-most visited site at time $t$ by an $X$ particle,  it was proved that
a.s. $r_t/t\to v$, with $v>0$ and that $\epsilon^{1/2}
(r_{\epsilon^{-1}t}-\epsilon^{-1}vt)$ converges to
a Brownian motion with non-degenerate variance.
In \cite{ks}, for the case of symmetric random walks with $D_X=D_Y$ a shape theorem was proved in
arbitrary dimensions.
In particular, in dimension $d=1$
 it was proved that if the initial configuration
of all the particles is a product Poisson measure with a finite
number of $X$ particles,  $r_t$ satisfies  a strong
law of large numbers.

Mai, Sokolov, Kuzovkov and Blumen \cite{mskb} performed Monte Carlo simulations
for a variation of the above described model in which both the $X$
and the $Y$ particles perform symmetric simple exclusion.
These numerical computations indicate that in the case $D_X=D_Y$,
with an initial condition which is a product Bernoulli measure, the front has
a ballistic movement with normal fluctuations.

In this paper we study a process
where the $X$ particles perform symmetric simple exclusion
but where  $D_Y=0$.
 We prove a law of large numbers and
a functional central limit theorem for the position of the
foremost visited site and for the number of activated
particles, giving an indication that corroborates
 the behavior observed in the numerical simulations
of \cite{mskb}.  In the model we consider, there are two
types of particles: the $X$ particles which move as
a symmetric simple exclusion process; the $Y$ particles, which do not move.
Initially there are no $X$ particles at sites $x>0$,
while the configuration $\eta:=\{\eta(0,x):x\le 0\}$ of
$X$ particles at sites $x\le 0$, where  $\eta(0,x)$  is the number of $X$ particles at site $x$,
is such that $\eta(0,0)=1$ but otherwise arbitrary.
Initially there are no $Y$ particles at $0,-1,-2,\ldots$ while
at sites $1,2,\ldots$, the $Y$ particles are distributed
according to 
 a product Bernoulli distribution of parameter $\rho$. 
When an $X$ particle jumps to a site where there is a $Y$ particle,
their position is interchanged and the $Y$ particle becomes
an $X$ particle. Since the $Y$ particles do not move, the dynamics
can be defined in terms of the configuration of the $X$ particles
and the rightmost visited site at time $t\ge 0$, which we call $r_t$. We
adopt the convention that $r_0=0$.
The state space of the process is then

$$
\Omega:=\left\{(r,p,\eta):r\in
\ze, p\in{\mathbb Z},\eta\in\{0,1\}^{\{\ldots,r-1,r\}}\right\}.
$$
Here $p$ represents a counter for the number of $Y$ particles which have been activated.
The infinitesimal generator of the process is

\begin{eqnarray}
\nonumber
& Lf(r,p,\eta):=\sum_{x,y\le r, |x-y|=1}\eta (x)(1-\eta(y))(f(r,p,
\sigma_{x,y}\eta)-f(r,p,\eta))\\
\nonumber
&+\rho\eta(r)(f(r+1,p+1,\eta+\delta_{r+1})-f(r,p,\eta))\\
\label{ex1}
&+(1-\rho)\eta(r)(f(r+1,p,\eta)-f(r,p,\eta)),
\end{eqnarray}
where $\delta_x$ denotes the configuration with one particle at
$x$, while $\sigma_{x,y}\eta$ denotes the configuration obtained
from $\eta$ after flipping the values of $\eta(x)$ and $\eta(y)$.

 We will use the
notation  $\eta(t,x)$ for the number of particles
at time $t\ge 0$ at site $x$, $p_t$ the value of the
counter for the  number
of activated $Y$ particles (so that $p_t-p_0$ is the actual number
of activated particles) and $r_t$ the position of the front
 at time $t\ge 0$. We will use the notation
 $\eta(t):=\{\eta(t,x):x\le r_t\}$ for the configuration of particles at time
$t$ and  will call the process $\{(r_t,p_t,\eta(t)):t\ge 0\}$
the {\it exclusion reactive process}.

\smallskip

\begin{theorem}
\label{theorem1} Assume that initially $r=p=0$, $\eta(0,0)=1$ and
$\eta(0,x)\in\{0,1\}$ is arbitrary for $x<0$.

\begin{itemize}

\item[(i)] There exist a $v>0$ and a $w>0$ which do not depend on the 
initial condition
$\{\eta(0,x): x\le 0\}$,  such that a.s.

$$
\lim_{t\to\infty}\frac{r_t}{t}=v,\qquad\qquad {\rm and}\qquad\qquad
\lim_{t\to\infty}\frac{p_t}{t}=w.
$$

\item[(ii)] There exist $\sigma_1$ and $\sigma_2$,
$0<\sigma_1, \sigma_2<\infty$, which do not depend on the initial condition
$\{\eta(0,x): x\le 0\}$, such that

$$\epsilon^{1/2}(r_{\epsilon^{-1}t}-\epsilon^{-1}vt),
\quad {\rm and}\quad
\epsilon^{1/2}(p_{\epsilon^{-1}t}-\epsilon^{-1}wt),\quad t\ge 0,$$
converge in law as $\epsilon\to 0$ to Brownian motions with
variances $\sigma_1$ and $\sigma_2$, respectively.

\end{itemize}
\end{theorem}

\begin{theorem}
\label{theorem2}     Consider the process as seen from the
front, $
\tau_{-r_t}\eta(t)
$.  There exist exactly two invariant measures: One 
 supported on the configuration
with no particles, and another, $\mu_\infty$.  The domain of attraction of the first consists of exactly the configuration with no 
particles.  Any nontrivial configuration in $\{0,1\}^{\{\ldots,-1,0\}}$ is in 
the domain of the second; if we denote by $\mu_t$ the distribution 
of the process $
\tau_{-r_t}\eta(t)
$, then $\mu_t\to \mu_\infty$ in the sense of weak convergence of
probability measures.  
\end{theorem}

\smallskip

As in \cite{cqr2}, 
  this model does not satisfy any obvious sub-additivity property which
would give a direct proof of part $(i)$ of theorem \ref{theorem1}.
But in a certain sense, the interaction given by the exclusion
dynamics of this paper, 
is stronger than the independent random walk dynamics of
\cite{cqr1}, or the annihilating random walks of \cite{cqr2}.

The proof of these results is based in regeneration time methods
as discussed in \cite{cqr2} and \cite{cqr1}. 
Nevertheless, to reduce the tail estimates of these regenerations times
to manageable expressions, we have 
 mapped the 
exclusion reactive process to a zero-range reactive
dynamics, with total jump rate $1$ per site.
 This mapping and the definition
of the regeneration times is given in section \ref{zerorange}.

In \cite{cqr2} and \cite{cqr1}, the regeneration times
where defined following \cite{sz}, in terms of two
alternating sequences of
stopping times. One of these sequences defines the first times
at which activated particles behind the leading one
branch. This definition is performed using a space-time line
which decouples the behavior of the old particles
behind the front with respect to the behavior of the front
itself. This is the approach  used in this paper
to define the regeneration times in section \ref{zerorange}.
However, we consider important to have an understanding
from a more fundamental point of view about alternative ways
which could be used to define regeneration times 
within the context of interacting particle systems
representing reactive fronts.
For this reason, in section \ref{alternative}, we present a definition
of regeneration times for the exclusion reactive process, which
is not done in terms of sequences of stopping times and which does
not require a space-time line decoupling the dynamics at the left
from the dynamics at the right of the front. 
The approach we present is in the spirit of Kesten \cite{k} within the context of 
Random Walks in Random Environments.

In section  \ref{p12}, theorems \ref{theorem1} and \ref{theorem2} are
proved assuming that the regeneration times and
the corresponding position of the front have finite
second moments. In section \ref{tail}, these estimates are
performed. The non-degeneracy of the variances $\sigma_1$ and $\sigma_2$
of Theorem 1 is proved in section \ref{nondegeneracy}.

\smallskip

%---------------------------------------------------------------------------------------------------------EXCLUSION-

%-------------------------------------------------------------------------------------------------------------------

\smallskip
\section{Zero-range reactive process and regeneration times}
\label{zerorange}

The first step in the proof of theorems \ref{theorem1} and
\ref{theorem2}, will be to 
couple the exclusion reactive process defined by 
(\ref{ex1}), with a zero range
process where particles branch at the right-most visited site.
Then, it will be enough to prove a law of large numbers and a
functional central limit theorem for the right-most visited site
of this new process, together with the convergence towards an invariant
measure for the law of this process as seen from its front. In the first subsection we will
define the zero-range reactive process. In subsection
\ref{lp}, we will construct a version of the zero-range
reactive process where particles are labeled. In subsection
\ref{auxiliary}, we construct an auxiliary process, which
gives lower bounds  for the position of the front. Then,
in subsection \ref{rt}, we will
define the regeneration times for the labeled version of
the zero-range reactive process following some
of the methods introduced in \cite{k}.

\subsection{One dimensional zero-range reactive process}
\label{odzrp}
 Consider a configuration
$(r,p,\eta)$ of the stochastic combustion process with exclusion.
Let $q:=r-p$. We define $\zeta(q)$ as the number of  particles between the site
$r+1$ and the rightmost empty site. In other words,
$\zeta(q):=r-x_1$, where $x_1:=\sup\{x\le 0:\eta(x)=0\}$.
Next define $x_2:=\sup\{x<x_1:\eta(x)=0\}$, the position
of the second right-most empty site, and $\zeta(q-1):=x_1-x_2-1$,
the number of particles between $x_1$ and $x_2$.
In general, for $n>2$, we define $x_n:=\sup\{x<x_{n-1}:\eta(x)=0\}$,
while $\zeta(q-n+1):=x_{n-1}-x_n-1$. Let
$\zeta(t):=\{\zeta(x,t):x\le q\}$ and $q_t:=r_t-p_t$ with
initial condition $(q_0,\zeta_0)=(q,\zeta)$. It is easy to check that the
stochastic process $\{(q_t,p_t,\zeta(t)), t\ge 0\}$ follows the dynamics of a one-dimensional
reactive process with a zero-range dynamics with infinitesimal
generator

\begin{eqnarray}
\nonumber
&Lf(q,p,\zeta)=\sum_{x,y\le q :|x-y|=1 } 1(\eta(x)>0)(f(q,p,\zeta^{x,y})-f(q,p,\zeta))\\
\nonumber
&+(1-\rho)1(\eta(q)>0)(f(q+1,p,\zeta-\delta_q+\delta_{q+1})-f(q,p,\zeta))\\
\label{coupling}
&+\rho 1(\eta(q)>0)(f(q,p+1,\zeta+\delta_q)-f(q,p,\zeta)).
\end{eqnarray}
This is a zero-range process with total jump rate $2$ at those
sites strictly to the left of the rightmost visited
site $q$,  with jump rate $1$ to the left and $1-\rho$ to the right at site $q$,
and with branching  at rate $\rho 1(\zeta(q)>0)$ at site $q$. 
Note that $p_t$ represents the number of times a branching has occurred.
We will call the triple $\{(q_t,p_t,\zeta(t)):t\ge 0\}$ the {\it
zero-range reactive process}.

\smallskip

\subsection{Labeled process}
\label{lp} We will make an
explicit construction of the zero-range reactive process where each
particle carries a label $z\in \ze$,  representing the priority
it has. The movement of a given particle will not be affected
by particles with smaller labels. 
The construction will be performed in
terms of a stochastic process $\{({\mathcal Y}(t),  q_t):t\ge 0
\}$, where at time $t\ge 0$, the first component ${\mathcal Y}(t)$ represents the
positions on $\ze$ of a random number of particles,  and $q_t$ the
position of
the rightmost visited site. 
Thus, the
state space of this process is ${\mathbb S}:={\mathbb T}\times  \ze$,
where ${\mathbb T}:=\cup_{A\subset\ze} \ze^A$. A typical
element of ${\mathbb T}$ will be denoted by $\{y_x:x\in A\}$, where
$A$ is the set of labels. We will furthermore use the notation
$p:=\sup\{x:x\in A\}$ and $p_t$, when the corresponding set $A_t$
is time dependent.

Let us fix an initial condition
$(\{y_x:x\in A_0\},q_0)\in{\mathbb S}$. 
Now associate to each  
$x\in A_0\cup\{z\in\ze: z>p_0\}$  a discrete time simple symmetric random walk
$X_x$ starting from $0$,  and a sequence $\{\tau_x^{(i)}:i\ge 1\}$ of
i.i.d. rate $2$ exponential random variables, 
which will represent the
potential jump times of an associated continuous time random walk.
Let us also choose another sequence $\{\upsilon_i:i\ge 1\}$ of i.i.d. 
Bernoulli random variables of parameter $\rho$. We choose all these random variables independent
of each other.

Let us first  define the dynamics of our process for an initial
condition $(\{y_x:x\in A_0\},q_0)$ such that $A_0\subset\ze$ is finite,
and such that $q_0=\sup\{y_x:x\in A_0\}$.
 We will
associate to each discrete time random walk $X_{x}, x\in A_0$,
a continuous time random walk $Y_{x}$, such that
$Y_x(0)=y_x$.
Let $n_1$ be the cardinality of the
 set $\{y_x:x\in A_0\}$.
 We identify at each site in this set, the particle with the
smallest label: let us call them $x_{1},\ldots ,x_{n_1}$. Let $x_*$ be the label within
the group $x_{1},\ldots,x_{n_1}$ where the minimum
 $\tau_1:=\min\{\tau_{x_{1}}^{(1)},\ldots,
\tau_{x_{n_1}}^{(1)}\}$ is achieved.
If $y_{x_*}<q_0$, we let the continuous time random walk  $Y_{x_*}$ jump at time
 $\tau_1$ according to $X_{x_*}$, while
the other random walks do not move. Thus, the
first change in the process $\{(\{Y_x(t):x\in A_t\},q_t):t\ge 0\}$
occurs at time $\tau_1$ when $Y_{x_*}(\tau_1)=y_{x_*}+X_{x_*}(1)$.
If $y_{x_*}=q_0$, and if the discrete time random walk $X_{x_*}$
jumps to the left at time $1$, we let $Y_{x_*}$ jump to the
left at time $\tau_1$.
Thus $Y_{x_*}(\tau_1)=y_{x_*}+X_{x_*}(1)$.
 If $y_{x_*}=q_0$, but
the  random walk $X_{x_*}$
jumps to the right at time $1$, we let $Y_{x_*}$ jump to the right
at time $\tau_1$ only if the Bernoulli random variable $\upsilon_1=0$,
case in which $Y_{x_*}(\tau_1)=y_{x_*}+X_{x_*}(1)$ and $q_{\tau_1}=q_0+1$, 
while if $\upsilon_1=1$, 
a random walk $Y_{p_0+1}$, which will follow the trajectory of
the  random walk $X_{p_0+1}+q_0$, is created
at time $\tau_1$ while the remaining random walks do not move, so that
$A_{\tau_1}=A_0\cup\{p_0+1\}$.
This defines
the dynamics of the process 
$\{({\mathcal Y}(t),q_t);t\ge 0\}$ in the time interval $[0,\tau_1]$,
with ${\mathcal Y}(t):=\{Y_x(t):x\in A_t\}$..

Let us now recursively define the process for arbitrary times.
Assume that for some $k$, such that $k\ge 1$, $n_{k}$ and $\tau_k$ have been defined
and also the process in the time interval $[0,\tau_{k}]$.
Let $p_{k}:=\sup\{x:x\in A_{\tau_{k}}\}$. Call $n_{k+1}$  the cardinality of the set 
$\{Y_x(\tau_{k}): x\in A_{\tau_{k}}\}$,  and identify at each site
in this set 
the particle with smallest label: we call these
labels $x_{1},\ldots, x_{n_{k+1}}$. Let $x_*$ be
the label where the minimum
 $\tau_{k+1}:=\min\{\tau_{x_{1}}^{(k+1)},\ldots,\tau_{x_{n_k}}^{(k+1)}\}$
is achieved. Denote for $1\le i\le n_k$ as $N_{x_i}$ the total
number of jumps performed up to time $\tau_k$ by the random walk $Y_{x_i}$.
If $Y_{x_*}(\tau_{k})<q_{\tau_{k}}$, we let the continuous time
 random walk $Y_{x_*}$ jump at time
 $\tau_{k+1}$ according to $X_{x_*}$ so that $Y_{x_*}(\tau_{k+1})
=y_{x_*}+X_{x_*}(N_{x_*}+1)$, while
the other random walks do not move.
If $Y_{x_*}(\tau_{k})=q_{\tau_{k}}$ and
if the random walk $X_{x_*}$ jumps to the left
at the time $N_{x_*}+1$, we let $Y_{x_*}$ jump
to the left at time $\tau_{k+1}$. 
If $Y_{x_*}(\tau_{k})=q_{\tau_{k}}$ but
the random walk $X_{x_*}$ jumps to the right
at the time $N_{x_*}+1$, we let $Y_{x_*}$ jump
to the right at time $\tau_{k+1}$ only if the Bernoulli
random variable $v_{1+p_{k}-p_0}=0$, while
if $v_{1+p_{k}-p_0}=1$
a random walk $Y_{p_k+1}$, following the discrete time trajectory of
$q_{\tau_{k}}+X_{p_k+1}$, is created.

Let us now consider the case in which $A_0$ 
 is not necessarily finite, but nevertheless
$p_0=\sup\{x:x\in A_0\}<\infty$.
Consider the initial condition of positions of the particles ${\mathcal Y}(0)=\{y_x:x\in A_0\}$
and of the front $q_0$.
We use the following notations:
given $A\subset\ze$ and $n\ge 1$, 
define  ${\mathcal Y}^n(0):=\{y_x\in{\mathcal Y}(0):y_x\le q_0-n\}$ and
$A^n:=\{x\in A: y_x\le q_0-n\}$.
Since $A^n_0$ is finite,
we can define the
process $\{({\mathcal Y}^n(t), q^n_t):t\ge 0\}$ as
in the previous paragraphs.

\smallskip
\begin{lemma}
\label{zrgd} There exists a set of full measure such that for every $t\ge 0$ the
following statements are true.

\begin{itemize}

\item[(i)]  There is an $n_0$ such that if $n\ge n_0$,

$$
A^n_t=A^{n_0}_t.
$$

\item[(ii)] Let $p^n_t:=\sup\{x:x\in A^n_t\}$. Then
$$
p_t:=\lim_{n\to\infty}p^n_t,
$$
exists.

\item[(iii)] For every $x\in A_t:=A^{n_0}_t$.
$$
Y_x(t):=\lim_{n\to\infty} Y^n_x(t),
$$
exists.

\item[(iv)] The limit
$$
q_t:=\lim_{n\to\infty}q^n_t,
$$
exists.

\end{itemize}
\end{lemma}

%CHECK THIS PROOF!!

\begin{proof} 
Let us prove part $(i)$. 
Without loss of generality we assume that for each $n\ge 1$,
$A_0^{n+1}\ne A_0^n$.
Consider the event $E_n:=\{A_t^m\ne A_t^n:{\rm for}\ {\rm some}\ m\ge n \}$.
By the lemma of Borel-Cantelli it is enough to prove that

\begin{equation}
\label{en}
\sum_{n=1}^\infty P[E_n]<\infty.
\end{equation}
Now the probability of the event $E_n$ is upper bounded by the
probability that in the time interval $[0,t]$, a Poisson process
of rate $2$ has performed at least $n$ steps. Indeed, the event $E_n$ is
contained in the event that some particle initially at a distance
larger than 
$n$ from the front, is alone at the foremost visited site before time $t$. But this can happen
only if at least $n$ jumps where performed 
before time $t$. Thus,

$$
P[E_n]\le \sum_{k=n}^\infty e^{-2t}\frac{(2t)^k}{k!}\le\frac{1}{n!},
$$
which proves (\ref{en}). Similar arguments can be used to prove parts $(ii)$, $(iii)$ and $(iv)$.
\end{proof}

\smallskip

Define now ${\mathcal Y}(t):=\{(Y_x(t): x\in A_t\}$. Then 
 the triple $\{({\mathcal Y}(t),q_t):t\ge 0\}$, defines a
probability measure ${\mathbb P}$ on the Skorokhod
space $D([0,\infty);{\mathbb S})$, which we will call the
{\it labeled zero-range reactive process}. It turns out that
the particle count

$$
\zeta(t,x)=\sum_{x'\in A_t} 1(Y_{x'}(t)=x),
$$
together with the pair $p_t$ and $q_t$ defined in lemma \ref{zrgd},
satisfies the dynamics defined by the infinitesimal generator
(\ref{coupling}).  

Consider an initial condition $w:=(\{y_x:x\in A_0\},q_0)\in{\mathbb S}$,
such that $y_{p_0}=q_0$, where $p_0=\sup\{x:x\in A\}$. Define $\delta_{p_0,q_0}
\in{\mathbb S}$  as $\delta_{p_0,q_0}:=(\{y_{x}:x\in  B_0\},q_0)$,
where $ B_0:=\{p_0\}$. This defines
two coupled labeled zero-range reactive processes
$\{(\{Y_x(t):x\in A_t\},q_t):t\ge 0\}$
and $\{(\{Y'_x(t):x\in  B_t\},q'_t):t\ge 0\}$, with initial conditions
$w$ and $\delta_{p_0,q_0}$ respectively.
The corresponding particle counts define two
coupled  zero range reactive processes with initial conditions
corresponding to $w$ and $\delta_{p_0,q_0}$, which we will
denote by $\zeta_w$ and $\zeta_{\delta_{p_0,q_0}}$ respectively.
The corresponding rightmost visited sites of these
processes will be denoted by $q^w_t$ and $q^{\delta_{p_0,q_0}}_t$,
whereas the counters of the number of activated particles by $p^w_t$ and $p^{\delta_{p_0,q_0}}_t$.
Furthermore we define

$$
\xi(t,x):=\sum_{x'\in  A_0-\{p_0\}} 1(Y_{x'}(t)=x)
$$
and

$$
\zeta'_w(t,x):=\sum_{x'\in  A'_t} 1(Y_{x'}(t)=x),
$$
where $A'_t:=(A_t-A_0)\cup\{p_0\}$.
The processes $\{(\zeta_w(t),q^w_t):t\ge 0\}$,
$\{(\zeta_{\delta_{p_0,q_0}}(t), q^{\delta_{p_0,q_0}}_t):t\ge 0\}$,
$\{\zeta'_w(t):t\ge 0\}$, and 
$\{\xi(t):t\ge 0\}$  are then coupled.

Let $V$ be the first time that some particle with a label smaller than $p_0$
is at the right-most visited site while no particle with label larger than
or equal to $p_0$ is at the right-most visited site

\begin{equation}
\nonumber
V':=\inf\left\{t\ge 0:\xi(t,q_t)> 0, \zeta_{\delta_{p_0,q_0}}(t,q_t)=0
\right\}.
\end{equation}

\smallskip

\begin{lemma} 
\label{lemma2}
 Let $A\subset {\bf Z}$,
$w:=(\{y_x:x\in A_0\},q_0)\in{\mathbb S}$
and $p_0=\sup\{x: x\in A_0\}$. Assume that
$y_{p_0}=q_0$. Consider the corresponding coupled process
$\{((\zeta_w(t),q^w_t),(\zeta_{\delta_{p_0,q_0}}(t), 
q^{\delta_{p_0,q_0}}_t),\xi(t));t\ge 0\}$.
Then

$$ q^{\delta_{p_0,q_0}}_t= q^w_t,\qquad t< V',$$ 
$$ p^{\delta_{p_0,q_0}}_t= p^w_t,\qquad t< V'$$ 
and

$$
\zeta_{\delta_{p_0,q_0}}(t)
=\zeta'_w(t), \qquad t< V'.
$$
\end{lemma}
\begin{proof} It is enough to observe that before time $V'$, none
of the particles with labels smaller than $p_0$ affect the
dynamics of the process $\{(\{Y_x(t):x\in A_t\}, q_t):t\ge 0\}$.
\end{proof}

\smallskip

\subsection{Auxiliary process}
\label{auxiliary} Let us now define a process
which will be helpful to obtain estimates for the law
of some stopping times used to define the regeneration times.
Let $n\ge 2$ be a fixed natural number. Now let $\{Z_x:x\ge 0\}$
be a set of independent continuous times simple symmetric
random walks of rate $2/n$, such that $Z_x(0)=x$. Define $\nu_1$ as the first time
the random walk $Z_0$ hits the site $1$,

$$
\nu_1:=\inf\{t\ge 0: Z_0(t)=1\}.
$$ 
For $m$ such that $2\le m\le n$, define $\nu_m$ as the first
time that any of the random walks $Z_0, Z_1,\ldots, Z_{m-1}$ hits
site $m$,

$$
\nu_m:=\inf\{t\ge 0: \sup_{0\le i\le m-1}Z_i(t)=m\}.
$$ 
And for $m>n$, define

$$
\nu_m:=\inf\{t\ge 0: \sup_{m-n\le i\le m-1}Z_i(t)=m\}.
$$
Now define for $t\ge 0$,

$$
\tilde q_t:=m\qquad {\rm for}\qquad 
\sum_{i=0}^{m  }\nu_i\le t< \sum_{i=0}^{m+1}\nu_i,
$$ 
with the convention that $\nu_0=0$. In the sequel we will
call $\{\tilde q_t: t\ge 0\}$, the {\it auxiliary process} with $n$
particles.

\smallskip
\begin{lemma} Whenever $n\ge 3$, there is an $\alpha>0$ such that

$$
\liminf_{t\to\infty}\frac{ q^{\delta_{0,0}}_t}{t}\ge \alpha\qquad a.s.
$$
\end{lemma}
\begin{proof} For each natural $n$, we can construct a process having the same law as 
$\{q_t^{\delta_{0,0}}:t\ge 0\}$: the  last $n$ activated particles
have priority over the ones activated previously; nevertheless, if at a given
time there are $m$ particles from this group (of $n$ particles), each
one jumps at a rate $2/m$. We can couple this construction with
the auxiliary process with $n$ particles $\{\tilde q_t:t\ge 0\}$ in such a way
that $q_t^{\delta_{0,0}}\ge \tilde q_t$ (for the details of such a coupling
within a similar context see \cite{cqr2}). Now, it is easy to check
that

$$
\lim_{t\to\infty}\frac{ \tilde q_t}{t}=: \alpha\qquad a.s.
$$
\end{proof}

\smallskip

\subsection{Regeneration times}
\label{rt} Let us consider an initial condition
$(\{y_x:x\in A_0\},q_0)$ such that $y_{p_0}=q_0$. Let also $\alpha_1$ and $\alpha_2$ 
be such that $0<2\alpha_1<\alpha_2<\alpha$.
Define

$$
T:=\inf\{t\ge 0: Z_{p_0,q_0}(t)=q_0-1\},
$$
which is the first time that the leading particle at $q_0$ jumps
backwards. Now let

$$
\bar q^{\delta_{p_0,q_0}}_t=Z_{p_0,q_0}(t)\qquad{\rm for}\ 0\le t<T
$$
and $\bar q^{\delta_{p_0,q_0}}_t=q^{\delta_{p_0,q_0}}_t$ for $t\ge T$.
Define

$$
U:=\inf\{t\ge 0:\bar q^{\delta_{p_0,q_0}}_t-q_0< \lfloor\alpha_2 t\rfloor
\ {\rm or}\ \zeta_{\delta_{p_0,q_0}}(t,0)=0
\},
$$

$$
V:=\inf\{t\ge 0:\xi(t,\lfloor\alpha_1 t\rfloor+q_0)>0\}
$$
and

$$
D:=\min\{U,V\}.
$$
The definition of $U$ in terms of $\bar q^{\delta_{p_0,q_0}}_t$
instead of $q^{\delta_{p_0,q_0}}_t$ is necessary because we have
to avoid the possibility that before time $D$ some particle
in the configuration $\xi$ branches before the front
$q^{\delta_{p_0,q_0}}_t$ has moved. Furthermore, the
fact that $2\alpha_1<\alpha_2$ guarantees that 
for times $t>1/\alpha_1$ the function $\lfloor\alpha_2 t\rfloor$
is never equal to $\lfloor\alpha_1 t\rfloor$.
 Note that $U$, $V$ and $D$
are stopping times with respect to the natural
filtration $\{{\mathcal F}_t:t\ge 0\}$ of the labeled zero-range reactive
process.

Let us also define the first time $U$ and $V$ happen after
time $s\ge 0$,

$$
U\circ\theta_s:=\inf\{t\ge 0: \bar q_t^{\delta_{p_s,q_s}}-q_s< \lfloor\alpha_2 t\rfloor\},
$$

$$
V\circ\theta_s:=\inf\left\{t\ge 0:\xi_{w_s}(t,\lfloor\alpha_1 t\rfloor+q_s)> 0\}
\right\},
$$
and $D\circ\theta_s:=\min\{U\circ\theta_s,V\circ\theta_s\}$,
where $w_s:=(\{Y_x(s):x\in A_s\},q_s)$.

For each $y\in {\mathbb N}$, define the ${\mathcal F}_t$-stopping time

$$
N_y:=\inf\{t\ge 0: p_t=p_0+y\}.
$$
We now define sequences of stopping times $\{S_k:k\ge 0\}$
and $\{D_k:k\ge 1\}$ as follows. First let $S_0:=0$ and $R_0:=0$. 
Then define for $k\ge 0$,

$$
S_{k+1}:=N_{R_k+1},\qquad D_{k+1}:=D\circ\theta_{S_{k+1}}+S_{k+1},
\qquad R_{k+1}:=p_{D_{k+1}},
$$
where we adopt the convention that $p_{D_k}=\infty$ and
 $S_{k+1}=\infty$ in the event $D_k=\infty$. We similarly
define $U_k:=U\circ\theta_{S_k}+S_k$ and
$V_k:=V\circ\theta_{S_k}+S_k$ for $k\ge 1$. Let now

$$
K:=\inf\{k\ge 1:S_k<\infty, D_k=\infty\},
$$
and define the {\it regeneration time},

$$
\kappa:=S_K.
$$
As in \cite{cqr2}, $\kappa$ is not a stopping time. Define ${\mathcal G}$,
the information up to time $\kappa$, as the completion of the $\sigma$-algebra
 generated by events of the form $\{t\le\kappa\}\cap A$, with $A\in
{\mathcal F}_t$.

The following proposition will be proved in section \ref{tail}. We call an 
element $w\in{\mathbb S}$ {\it nontrivial}, if it corresponds to
a configuration with at least one particle behind the front or at it.

\smallskip

\begin{proposition}
\label{base}
 For every non-trivial initial condition $w\in{\mathbb S}$,

\begin{equation}
\label{kappa1}
\kappa<\infty,\qquad\qquad {\mathbb P}_w-a.s.
\end{equation}
Furthermore,

\begin{equation}
\label{kappa11}
{\mathbb E}_{\delta_{0,0}}[\kappa^2|U=\infty]<\infty\qquad {\rm and}\qquad
{\mathbb E}_{\delta_{0,0}}[r_\kappa^2|U=\infty]<\infty.
\end{equation}
\end{proposition}
\smallskip

\smallskip

\begin{proposition}
\label{prop1}
 Let $F$ be a Borel subset of $D([0,\infty);\Omega)$.
Then, for every nontrivial $w\in{\mathbb S}$,

$$
{\mathbb P}_w[(q_{\kappa+\cdot}-q_\kappa,
p_{\kappa+\cdot}-p_\kappa,\tau_{-q_\kappa,-p_\kappa}\zeta'_{w_\kappa}(\kappa+\cdot))
\in F|{\mathcal G}]
={\mathbb P}_{\delta_{0,0}}[
(q_\cdot,p_\cdot,
\zeta(\cdot))\in F|U=\infty].
$$
\end{proposition}

\begin{proof} It is enough  to prove that for every $B\in{\mathcal G}$,

\begin{eqnarray*}
&{\mathbb P}_w[B,(q_{\kappa+\cdot}-q_\kappa,
p_{\kappa+\cdot}-p_\kappa,\tau_{-q_\kappa,-p_\kappa}\zeta'_{w_\kappa}(\kappa+\cdot))
\in F]\\
&={\mathbb P}[B]
{\mathbb P}_{\delta_{0,0}}[
(q_\cdot,p_\cdot,
\zeta(\cdot))\in F|U=\infty].
\end{eqnarray*}
As in \cite{cqr2}, this can be done using  lemma \ref{lemma2} and
observing that on the event
on the event $S_k<\infty$, $q_{s_k}=x$ and $p_{s_k}=y$,
we have that

$$
\zeta'_{w_\kappa}(S_k+\cdot)=\zeta_{\delta_{x,y}}(\cdot),
$$
whenever $U_k=V_k=\infty$. 
\end{proof}

\smallskip

We can now define a sequence $\kappa_1\le\kappa_2\le\cdots$,
with $\kappa_1:=\kappa$ while for $n\ge 1$,

$$
\kappa_{n+1}:=\kappa_n+\kappa(w_{\kappa_n}+\cdot),
$$
where $\kappa(w_{\kappa_n}+\cdot)$ is the regeneration time
starting from $w_{\kappa_n}$ and we set $\kappa_{n+1}=\infty$
on the event $\kappa_n=\infty$. We call $\kappa_1$ the {\it 
first regeneration time}
and $\kappa_n$ the {\it $n$-th regeneration time.}
Now define for each $n\ge 1$, the $\sigma$-algebra ${\mathcal G}_n$
as the completion with respect to ${\mathbb P}$
of the smallest $\sigma$-algebra containing all sets of the
form
$\{\kappa_1\le t_1\}\cap\cdots\cap\{\kappa_n\le t_n\}\cap A$.
$A\in{\mathcal F}_{t_n}$. As in \cite{cqr2}, we then have the
following generalization of proposition \ref{prop1}.

\smallskip

\begin{proposition}
\label{prop2}
 Let $F$ be a Borel subset of $D([0,\infty);\Omega)$.
Then, for every nontrivial $w\in{\mathbb S}$ and 
$n\ge 1$,

\begin{eqnarray*}
&{\mathbb P}_w[(q_{\kappa_n+\cdot}-q_{\kappa_n},
p_{\kappa_n+\cdot}-p_{\kappa_n},\tau_{-q_{\kappa_n},-p_{\kappa_n}}
\zeta'_{w_{\kappa_n}}({\kappa_n}+\cdot))
\in F|{\mathcal G}_n]\\
&={\mathbb P}_{\delta_{0,0}}[
(q_\cdot,p_\cdot,
\zeta(\cdot))\in F|U=\infty].
\end{eqnarray*}
\end{proposition}
\smallskip

\begin{proposition}
\label{iid}
 Let $w\in {\mathbb S}$.
(i)  Under ${\mathbb P}_w$, $\kappa_1,\kappa_2-\kappa_1, \kappa_3-\kappa_2, \ldots$ are independent, and $\kappa_2-\kappa_1, \kappa_3-\kappa_2, \ldots$ are identically distributed with law identical to that of $\kappa_1$ under
${\mathbb P}_{\delta_{0,0}}[\cdot |U=\infty]$.
(ii) Under ${\mathbb P}_w$, $ q_{\cdot\land\kappa_1}, q_{(\kappa_1+\cdot)\land\kappa_2} - q_{\kappa_1}, q_{(\kappa_2+\cdot)\land\kappa_3} - q_{\kappa_2}, \ldots$ are independent, and $q_{(\kappa_1+\cdot)\land\kappa_2} - q_{\kappa_1}, q_{(\kappa_2+\cdot)\land\kappa_3} - q_{\kappa_2}, \ldots$ are identically distributed with law identical to that of $q_{\kappa_1}$ under
${\mathbb P}_{\delta_{0,0}}[\cdot|U=\infty]$ .
(iii) Under ${\mathbb P}_w$, $ p_{\cdot\land\kappa_1}, p_{(\kappa_1+\cdot)\land\kappa_2} - p_{\kappa_1}, p_{(\kappa_2+\cdot)\land\kappa_3} - p_{\kappa_2}, \ldots$ are independent, and $p_{(\kappa_1+\cdot)\land\kappa_2} - p_{\kappa_1}, p_{(\kappa_2+\cdot)\land\kappa_3} - p_{\kappa_2}, \ldots$ are identically distributed with law identical to that of $p_{\kappa_1}$ under
${\mathbb P}_{\delta_{0,0}}[\cdot|U=\infty]$ .
\end{proposition}

\smallskip

\section{Proof of theorems 1 and 2}
\label{p12}

\subsection{Proof of theorem 1}  The proof of parts $(i)$ and
$(ii)$ of theorem 1 follow now using standard arguments (see
for example  \cite{cqr2} or in the context of Random Walks in 
Random Environment 
\cite{sz}). Indeed, using Propositions  \ref{base} and \ref{iid}
we first prove that, a.s.

$$
\lim_{n\to\infty}\frac{q_{\kappa_n}}{\kappa_n}=
\frac{{\mathbb E}_w[q_{\kappa_1}|U=\infty]}
{{\mathbb E}_w[\kappa_1|U=\infty]}=:v_1
\qquad {\rm and}\qquad
\lim_{n\to\infty}\frac{p_{\kappa_n}}{\kappa_n}=\frac{{\mathbb E}_w[p_{\kappa_1}|U=\infty]}
{{\mathbb E}_w[\kappa_1|U=\infty]}=:v_2.
$$
An interpolation argument then proves that a.s.
$$
\lim_{t\to\infty}\frac{q_t}{t}=v_1
\qquad {\rm and}\qquad
\lim_{t\to\infty}\frac{p_t}{t}=v_2.
$$
Since $r_t=q_t+p_t$, this implies that a.s.

$$
\lim_{t\to\infty}\frac{r_t}{t}=v_1+v_2,
$$
proving part $(i)$ of theorem \ref{theorem1}.
We then define 
$P_j:=p_{\kappa_{j+1}}-p_{\kappa_j}-(\kappa_{j+1}-
\kappa_j)v_2$ and
$R_j:=p_{\kappa_{j+1}}+q_{\kappa_{j+1}}-p_{\kappa_j}-q_{\kappa_j}-(\kappa_{j+1}-
\kappa_j)(v_1+v_2)$,
and show that

$$
\Sigma_m:=\sum_{j=1}^m P_j\qquad {\rm and}\qquad \Sigma'_m:=\sum_{j=1}^m R_j,
$$
converge in law to Brownian motions with variances 

$$
\sigma_1:=\frac{{\mathbb E}_w[(p_{\kappa_1}-v_2\kappa_1)^2|U=\infty]}
{{\mathbb E}_w[\kappa_1|U=\infty]}
\qquad {\rm and}\qquad
\sigma_2:=\frac{{\mathbb E}_w[(r_{\kappa_1}-(v_1+v_2)\kappa_1)^2|U=\infty]}
{{\mathbb E}_w[\kappa_1|U=\infty]}
$$
respectively.
An interpolation argument can then be used to obtain the full
limits proving that

\begin{equation}
\label{bm}
\epsilon^{1/2}(p_{\epsilon^{-1}t}-\epsilon^{-1}t v_2)\qquad
{\rm and}\qquad
\epsilon^{1/2}(r_{\epsilon^{-1}t}-\epsilon^{-1}t (v_1+v_2)),\quad t\ge 0,
\end{equation}
converge in law to Brownian motions with variances $\sigma_1$ and
$\sigma_2$ respectively. In section \ref{nondegeneracy} we will prove
that these variances are positive.
\smallskip

\subsection{Proof of theorem 2}  Let $\mu_t$ be the law at time $t$
of the exclusion reactive process seen from the front
$\tau_{-r_t}\eta_t\in\Omega_0:=\{0,1\}^{{\mathbb Z}_-}$.
This is a Markov process with infinitesimal
generator

\begin{eqnarray*}
&L_0 f(\eta )=\rho\eta(0)(f(\tau_{-1}\eta+\delta_0)-f(\eta ))
+
(1-\rho)\eta(0)(f(\tau_{-1}\eta)-f(\eta))\\
&
+\sum_{x,y\le 0, |x-y|=1}\eta(x)(1-\eta(y))(f(\sigma_{x,y}\eta)
-f(\eta)).
\end{eqnarray*}
 For
a local function $f$ on $\Omega_0$, define $l(f)$ as
the smallest integer $l$ such that $f(\eta)$ does not depend
on $\eta(x)$ if $x<-l$. Now define the probability measure
$\mu_\infty$ on $\Omega_0$ by the formula

$$
\int_{\Omega_0} f d\mu_\infty=
\frac{{\mathbb E}_{\delta_0}[\int_{\kappa_N}^{\kappa_{N+1}} 
f(\tau_{-r_s}\eta_s)ds|U=\infty] }{{\mathbb E}_{\delta_0}[\kappa_1|U=\infty]},
$$
for $N$ such that $N(\alpha'-\alpha'')>l(f)$. This defines a consistent family
of probability measures on cylinders.

\smallskip

\begin{theorem} $\lim_{t\to\infty}\mu_t=\mu_\infty$ weakly
and $\mu_\infty$ is invariant for the generator $L_0$.
\end{theorem}
\begin{proof} First note that for $f$ local and $w$ an arbitrary initial
condition,
$\lim_{t\to\infty}{\mathbb E}_w[\kappa_{N+1}>t, f(\tau_{-r_t}
\eta(t))]=0$. Therefore in the decomposition

$$
\int_{\Omega_0} fd\mu_t=
{\mathbb E}_w[\kappa_{N+1}\le t, f(\tau_{-r_t}
\eta(t))]+
{\mathbb E}_w[\kappa_{N+1}>t, f(\tau_{-r_t}
\eta(t))],
$$
it is enough to examine the first term. Now, as in \cite{cqr2},
by Proposition \ref{prop2}
\begin{equation}
{\mathbb E}_w[\kappa_{N+1}\le t, f(\tau_{-r_t}
\eta(t))]=\int_0^t{\mathcal N}_t(d u) F_f(u),
\end{equation}
where

$$
{\mathcal N}_t([0,u]):=\sum_{k\ge 1}{\mathbb P}_w (\kappa_k\in [t-u,t])
$$
and

$$
F_f(u):={\mathbb E}_{\delta_{0,0}}[\kappa_N\le u<\kappa_{N+1}, f(\tau_{-r_t}
\eta(u))|U=\infty].
$$
Also, as in \cite{cqr1}, we can check the spread out assumption of the renewal
theorem (theorem 6.2 of \cite{th}) to show that

$$
\lim_{t\to\infty}{\mathcal N}_t([0,u])=\frac{u}{{\mathbb E}_{\delta_{0,0}}[\kappa_1
|U=\infty]},
$$
uniformly on compacts. Finally, since $|F_f|\le ||f||_\infty|F_1|$
and $\int F_1 du<\infty$ we conclude that

$$
\lim_{t\to\infty}\int_{\Omega_0} f d\mu_t=\int_{\Omega_0} f d\mu_\infty.
$$

\end{proof}

\smallskip
\section{Estimates for the regeneration times}
\label{tail}

Let us first state the  following estimate for the stopping time $U$.

\smallskip

\begin{lemma}
\label{ut}
 For every $p\ge 1$ there exists a constant $C=C(p)$
such that for every initial condition $w$

$$
{\mathbb P}_w[t\le U<\infty]\le C t^{-p}.
$$
\end{lemma}

\begin{proof} The proof is based on a comparison with the auxiliary process
defined in subsection \ref{auxiliary}, and  large deviation
estimates as in the proof of lemma 7 of \cite{cqr2}.
\end{proof}

\smallskip

Let us define $\tilde V$ as the first time that
independent continuous time simple symmetric random
walks created at rate $1$ at the
origin, reach the site $\lfloor \alpha t\rfloor$.
We  have the following estimate for the stopping time $V$.

\smallskip

\begin{lemma}
\label{vt} There is a constant $C>0$ such that for
every initial condition $w$

$$
{\mathbb P}_w[t\le V<\infty]\le \exp\left\{-Ct\right\}.
$$
\end{lemma}
\begin{proof} Let us note using the labeled construction of the zero-range
reactive process, that we can bound
the probability of  event $\{t\le V<\infty\}$ by the 
probability of the event $\{t\le \tilde V<\infty\}$. It is not
difficult to check that there is a constant $C>0$ such that,
$P[t\le \tilde V<\infty]\le\exp\left\{-Ct\right\}$ (see \cite{cqr2}).
\end{proof}

\smallskip

From Lemmas \ref{ut} and \ref{vt} we can now
directly prove (see \cite{cqr2}),
\smallskip
\begin{lemma}
\label{dtail}
 For every $p\ge 1$, there is a constant $C=C(p)$ such that for
every initial condition $w$

$$
{\mathbb P}_w[t\le D<\infty]\le C t^{-p}.
$$
\end{lemma}
\smallskip

We continue with two important lemmas. The first one is proved
using the auxiliary process as in \cite{cqr2}.

\smallskip

\begin{lemma}
\label{uinf}
 There is a $\delta_1>0$ such that for every
initial condition $w$

$$
{\mathbb P}_w[U=\infty]\ge \delta_1.
$$
\end{lemma}

\smallskip
The second lemma which follows can be proved using the inequality
${\mathbb P}_w[V=\infty]\ge{\mathbb P}_w[\tilde V=\infty]$.

\smallskip

\begin{lemma}
\label{vinf}
 There is a $\delta_2>0$ such that for every
initial condition $w$

$$
{\mathbb P}_w[V=\infty]\ge \delta_2.
$$
\end{lemma}

\smallskip

We will also need the following lemma.

\smallskip

\begin{lemma}
\label{qp}
 There is a constant $C>0$ such that for every
initial condition $w$ with $p=q=0$ and $M>1$,

$$
{\mathbb P}_w[q_t\ge Mt]\le C\exp\left\{-Ct\right\}\qquad
{\rm and}\qquad
{\mathbb P}_w[p_t\ge Mt]\le C\exp\left\{-Ct\right\}.
$$
\end{lemma}
\begin{proof} It is enough to note that the rate at which the process
$\{q_t:t\ge 0\}$ increases is always bounded by $1$. Similarly
for the process $\{p_t:t\ge 0\}$.
\end{proof}

\smallskip

Let us now prove (\ref{kappa1}) of Proposition \ref{base}.
Note that for
every $k\ge 1$,  ${\mathbb P}_w[\kappa=\infty]\le {\mathbb P}_w[D_k<\infty]$.
Taking the limit when $k\to\infty$ we obtain (\ref{kappa1}).
Part (\ref{kappa11}) of Proposition \ref{base} is a
consequence of the following lemma.

\smallskip

\begin{lemma} There is a constant $C>0$ such that for
every $p\ge 1$ and $t\ge 0$,

$$
{\mathbb P}_{\delta_{0,0}}[\kappa>t|U=\infty]\le C t^{-p}.
$$
\end{lemma}

\begin{proof} Note that

$$
{\mathbb P}_{\delta_{0,0}}[\kappa>t|U=\infty]=\sum_{k=1}^\infty
{\mathbb P}_{\delta_{0,0}}[S_k>t, K=k|U=\infty].
$$
From Lemmas \ref{uinf} and \ref{vinf}, 
using the strong Markov property we deduce that
there is a $\delta>0$ such that for every $l\ge 1$,

\begin{equation}
\label{kp}
{\mathbb P}_{\delta_{0,0}}[\kappa>t|U=\infty]\le
\sum_{k=1}^l
{\mathbb P}_{\delta_{0,0}}[S_k>t, K=k|U=\infty]+
\delta^{-1} (1-\delta)^l,
\end{equation}
Let $0<\gamma<1$ and consider the event

$$
A_k:=\{r_{D_1}-r_{S_1}<t^\gamma,\ldots, r_{D_{k-1}}-r_{S_{k-1}}<t^\gamma\}.
$$
On $A_k$ we have $r_{S_k}\le k t^\gamma$. Since $\bar q_t=q_t$ (because
the initial condition is $\delta_{0,0}$), if $U=\infty$ then
$r_t\ge q_t>\lfloor \alpha t\rfloor$ for all $t\ge 0$. We then have on
$A_k\cap\{U=\infty\}$ that $\lfloor\alpha S_k\rfloor\le k t^\gamma$.
Hence, for $t>(lt^\gamma+1)/\alpha$ and $k\le l$,

$$
{\mathbb P}_{\delta_{0,0}}[t<S_k<\infty, A_k|U=\infty]=0
$$
and

\begin{eqnarray}
\nonumber
&{\mathbb P}_{\delta_{0,0}}[t<S_k<\infty|U=\infty]\le
{\mathbb P}_{\delta_{0,0}}[t<S_k<\infty, A^c_k|U=\infty]\\
\label{sk}
&\le C\sum_{i=1}^{k-1}{\mathbb P}_{\delta_{0,0}}[r_{D_i}
-r_{S_i}\ge t^\gamma, S_k<\infty],
\end{eqnarray}
where we have used Lemma \ref{uinf}. Now, using Lemmas \ref{qp} and
\ref{dtail},
we can show that for every $p\ge 1$ there exists a constant $C$
such that

$$
{\mathbb P}_{\delta_{0,0}}[r_{D_i}
-r_{S_i}\ge t^\gamma, S_k<\infty]\le C t^{-p}.
$$
Substituting this estimate back into (\ref{sk}) and then into (\ref{kp}),
we conclude the proof.
\end{proof}

\smallskip

\section{Non-degeneracy}
\label{nondegeneracy}

 To prove the non-degeneracy of the limit of
the first expression in display (\ref{bm}) we will prove that for some 
$\alpha<\beta<v_2$ it is true that

\begin{equation}
\label{pk0}
{\mathbb P}_{\delta_{0,0}}[p_{\kappa_1}=1, \beta^{-1}<\kappa_1|U=\infty]>0.
\end{equation}
Note that,

$$
{\mathbb P}_{\delta_{0,0}}[p_{\kappa_1}=1, \beta^{-1}<\kappa_1,U=\infty]\ge
{\mathbb P}_{\delta_{0,0}}[\beta^{-1}<S_1<U , D\circ\theta_{S_1}=\infty].
$$
The lower bound of the above inequality can be written as

$$
{\mathbb E}_{\delta_{0,0}}[1(\beta^{-1}<S_1<U)
{\mathbb E}_{\delta_{0,0}}[ U\circ\theta_{S_1}=\infty,
V\circ\theta_{S_1}=\infty|{\mathcal F}_{S_1}]].
$$
Now

\begin{eqnarray*}
&{\mathbb E}_{\delta_{0,0}}[ U\circ\theta_{S_1}=\infty,
V\circ\theta_{S_1}=\infty|{\mathcal F}_{S_1}]=\\
&{\mathbb E}_{\delta_{0,0}}[ U\circ\theta_{S_1}=\infty|{\mathcal F}_{S_1}]
{\mathbb E}_{\delta_{0,0}}[ 
V\circ\theta_{S_1}=\infty|{\mathcal F}_{S_1}]\ge\delta_1\delta_2,
\end{eqnarray*}
where in the inequality we have used Lemmas \ref{uinf} and \ref{vinf}. Thus,

$$
{\mathbb P}_{\delta_{0,0}}[\beta^{-1}<S_1<U , D\circ\theta_{S_1}=\infty]
\ge\delta_1\delta_2 {\mathbb P}_{\delta_{0,0}}[\beta^{-1}<S_1<U].
$$
Now, the probability of the right-hand term can be lower bounded
by the probability that the particle initially at site $0$ with
label $0$, performs a first jump at some time $t$ towards the right
and then at time $t+t'$ it branches at site $1$
 such that $\beta^{-1}<t+t'<\alpha^{-1}$. Clearly this probability
is positive.

To prove the non-degeneracy of the second
 expression in display (\ref{bm}) this time we will show that for some
$\alpha<\beta<v_1+v_2$ it is true that

\begin{equation}
\label{pk}
{\mathbb P}_{\delta_{0,0}}[r_{\kappa_1}=1, \beta^{-1}<\kappa_1|U=\infty]>0.
\end{equation}
But note that $\{r_{\kappa_1}=1\}=\{p_{\kappa_1}=1\}$. Therefore
inequality (\ref{pk}) follows from (\ref{pk0}).

\smallskip
\section{An alternative definition of regeneration times}
\label{alternative}

\subsection{Asymmetric simple random walk}
 Let us show through a simple
example how to prove the independence of the increments of the
regeneration times, without defining them through a sequence
of stopping times. What we present here is essentially contained
in the paper of Kesten \cite{k}.

Consider a discrete time asymmetric simple random walk $\{X_n:n\ge
0\}$, which at each site jumps to the right with probability $q$
and to the left with probability $p$, with $q>p$. Define the
random time

$$\kappa:=\min\{n\ge 0:\min_{k\ge n}X_k>\max_{j<n}X_j\},$$
with the convention that $\max_{j<0}X_j=-1$. In words, $\kappa$ is
the first time the random walk visits a new site, without never
afterward moving to the left of such a site. Let us define the
$\sigma$-algebra ${\mathcal F}_\kappa$ of $\kappa$ as the one
generated by events of the form $\{\kappa
=n\}\cap\{X_1=x_1,\ldots,X_n=x_n\}$ for $n\ge 0$ and $x_1,\ldots
,x_n$ in ${\mathbb Z}^d$. If we call $P$ the law of the random
walk $\{X_n\}$,  starting from $0$, the following equality is
satisfied for every subset $A\in {\mathbb Z}^{\mathbb N}$:

\begin{lemma}
\label{lemma1}
\begin{equation}
\label{srw} P\left[ X_{\kappa+\cdot}-X_\kappa \in A|{\mathcal
F}_\kappa\right]=P\left[X_\cdot\in A|\kappa =0\right].
\end{equation}
\end{lemma}
Let us prove this equality using directly the Markov property. For
convenience, we define $Y_n:=\sup_{k\le n} X_k$, for $n\ge 0$.
Note first that

\begin{eqnarray}\label{decomposition}
  \{\kappa =n\}=\{\inf_{m\ge n} X_m\ge X_n\}\cap\{\sup_{1\le p\le
  n-1}X_p=X_n-1\}\cap B_n,
\end{eqnarray}
where

\begin{eqnarray*}
B_n&:=&\left\{\forall k<n:\ {\rm there}\ {\rm exists}\ {\rm a}\
j,\ {\rm
such}\ {\rm that}\ k<j<n\ {\rm and}\ X_j<Y_k\right.\\
&\ &\left. {\rm or}\ {\rm there}\ {\rm exists}\ {\rm a}\ j,\ {\rm
such}\ {\rm that}\ 0\le j<k\ {\rm and}\ Y_j=X_k\right\}.
\end{eqnarray*}
Let us remark that $B_n\in {\mathcal F}_n$, where $\{{\mathcal
F}_n:n\ge 0\}$ is the natural filtration of the random walk $\{
X_n:n \ge 0\}$. Furthermore, the condition $q>p$, ensures that
$P[\kappa =n]>0$. Defining $C_n:=\{\sup_{1\le p\le n-1}X_p=X_n-1\}
\cap B_n\in {\mathcal F}_n$, we now have

\begin{eqnarray}
\nonumber &\ &P\left[ X_{\kappa+\cdot}-X_\kappa \in
A|\kappa=n,X_1=x_1,\ldots,
X_n=x_n\right]\\
&=& \nonumber \frac{P\left[ X_{\kappa+\cdot}-X_\kappa \in
A,\inf_{m\ge n}X_m\ge X_n, C_n,X_1=x_1,\ldots, X_n=x_n\right]}
{P\left[\inf_{m\ge n}X_m\ge
X_n,C_n,X_1=x_1,\ldots,X_n=x_n\right]}.
\end{eqnarray}
On the other hand, by the fact that $C_n\in{\mathcal F}_n$, the
Markov property, and translation invariance, we have that

\begin{eqnarray}
\nonumber &\ & P\left[ X_{\kappa+\cdot}-X_\kappa \in A,\inf_{m\ge
n}X_m\ge X_n,
C_n,X_1=x_1,\ldots, X_n=x_n\right]\\
\nonumber &=& E\left[{\bf
1}_{\{C_n,X_1=x_1,\ldots,X_n=x_n\}}P\left[\left.
X_{n+\cdot}-X_n\in A,\inf_{m\ge n} X_m\ge X_n\right|{\mathcal
F_n}\right]
\right]\\
\label{srw1} &=& P[C_n,X_1=x_1,\ldots,X_n=x_n]P[X_\cdot\in
A,\inf_{m\ge 0}X_m\ge 0].
\end{eqnarray}
Choosing $A$ as the whole space in the previous development, we
conclude that

\begin{eqnarray}
\nonumber
&\ & P\left[\inf_{m\ge n}X_m\ge X_n, C_n,X_1=x_1,\ldots, X_n=x_n\right]\\
\nonumber &=&P[C_n,X_1=x_1,\ldots,X_n=x_n]P[\inf_{m\ge n}X_m\ge
0],
\end{eqnarray}
which combined with (\ref{srw1}) completes the proof of
(\ref{srw}).

\smallskip

\subsection{Exclusion reactive process} Let us now show how the
approach presented in the context of a symmetric simple random walk
can be implemented for the exclusion reactive process to
define a version of the regeneration times.
To define the regeneration times we will consider the holes (empty sites) of the process as
second class particles.

We will construct the exclusion reactive process associating to
each bond connecting two nearest neighbor sites independent Poisson
processes each one of rate $1$. 
Let us assume that the initial condition is of the form
$(0,\eta(0))$, where $\eta(0)$ is any nontrivial configuration
of particles.  
If at a given time the rightmost visited site is $r$ and $x+1\le r$,
whenever the Poisson clock connecting sites $x$ and $x+1$ rings,
the state at sites $x$ and $x+1$ is interchanged: if the state
of the process was $\eta$, this is changed to $\sigma_{x,x+1}\eta$.
In this case the front stays at $r$.
If at a given time the rightmost visited site $r$ is occupied and
the Poisson clock connecting sites $r$ and $r+1$ rings,
with probability $\rho$ a particle is created at site $r+1$
and with probability $1-\rho$ the particle at $r$ jumps
to site $r+1$. In both cases the front advances one step to $r+1$.
Define $\tau_1$ as the first time the front
advances one step creating a new particle by

$$
\tau_1:=\inf\{t\ge 0: r_t>0\ {\rm and}\ \eta(r_t-1,t)=\eta(r_t,t)=1\}.
$$
In general for $n\ge 2$ define recursively $\tau_n$ as
the first time after time $\tau_{n-1}$ that the front advances
one step creating a new particle by

$$
\tau_n:=\inf\{t\ge \tau_{n-1}: r_t>r_{\tau_{n-1}}\ {\rm and}\
\eta(r_t-1,t)=
 \eta(r_t,t)=1\}.
$$
For each $n\ge 1$, at the time $\tau_n$ we will fill up the holes
defined by the configuration $\eta(\tau_n)$ with particles: the
original particles in this configuration are first class with
respect to the new ones, which we will call {\it holes}.
After time $\tau_n$ holes can never activate a new particle.

Let us now adopt the convention that whenever the Poisson clock of a bond 
corresponding to two occupied sites rings, the corresponding
particles or holes are interchanged. We can then label each particle and
each hole, following
its trajectory. We will assign the label $0$ to the particles
initially at $r$ or to the left of $r$. We assign the label $1$ to
the particle activated at time $\tau_1$. In general for $n\ge 2$,
we assign the label $n$ to the particle activated at time $\tau_n$.
Similarly, we assign the label $0$ to all the holes initially at
$r$ or to the left of $r$, and to the new holes created 
before time $\tau_1$.
 We assign the label $1$ to the  holes
created after time $\tau_1$ but before time $\tau_2$, and for $n\ge 2$ we
assign the label $n$ to the holes created after time $\tau_n$
but before time $\tau_{n+1}$.
Let us call $\Upsilon$ the state space
of this process consisting of the ordered pairs $(r,{\boldsymbol\eta})$,
where $r$ is an integer representing the position of the front and
${\boldsymbol\eta}$ is a configuration of labeled particles
and holes at sites $\ldots,r-1,r$. Let us call $Q$ the
law of this process in the corresponding Skorohod space.
We now define for each $n\ge 1$, the stopping time
$D_n$ as the first time after time $\tau_n$ that some of the
particles or holes with labels strictly smaller than $n$ (thus,
excluding the  foremost particle, created at time $\tau_n$)
is at the front

$$
D_n:=\inf\{t\ge 0:\tilde \eta(\tau_n+t,r_{\tau_n+t})=1\}.
$$
Here, at time $t\ge 0$, $\tilde\eta(\tau_n+t)$ is the particle-hole count
of those particles or holes in ${\boldsymbol\eta}(\tau_n+t)$ 
with labels strictly smaller than $n$. 
Now define  the first regeneration time
$$
\kappa_1:=\inf\{\tau_m: m\ge 1, D_m=\infty\}.
$$
We then define recursively for $n\ge 2$,

$$
\kappa_n:=\inf\{\tau_m: \tau_m>\kappa_{n-1}, D_m=\infty\}.
$$
Let us call $\delta_0\in\Upsilon$ any initial condition with
rightmost visited site $r=0$, one particle with label
$0$ at $0$ and none elsewhere, and one hole at each site
$x<0$ each one with label $-1$. Given this initial condition,
we define $D_0$ as the first time that one of the holes with label
$-1$ is at the front.
Call ${\mathcal G}_1$ the information up to time $\kappa_1$.
We have the following proposition corresponding to
Proposition \ref{prop1}. We will in general call $\eta$ the
particle count corresponding to a state $\boldsymbol\eta$
and
$\Omega':=\{(r,\eta):r\in{\mathbb Z},\eta\in\{0,1\}^{\{\ldots,r-1,r\}}\}.$

\smallskip

\begin{proposition}
 Let $F$ be a Borel subset of $D([0,\infty);\Omega')$.
Then, for every nontrivial ${\boldsymbol \eta}\in\Upsilon$,
which has only particles with label $0$ and holes with label $0$,
$$
Q_{\boldsymbol \eta}[(r_{\kappa_1+\cdot}-r_{\kappa_1},\tau_{-r_{\kappa_1}}\eta'({\kappa_1}+\cdot))
\in F|{\mathcal G}_1]
=Q_{\delta_{0}}[
(r_\cdot,\eta(\cdot))\in F|D_0=\infty],
$$
where $\eta'$ is the particle count of those particles
created at time $\kappa_1$ or after.
\end{proposition}
\begin{proof} Note that the events $A\cap\{\kappa_1=\tau_n\}$ with
$A\in{\mathcal F}_{\tau_n}$ for some $n\ge 1$,
generate the $\sigma$-algebra ${\mathcal G}_1$.
Therefore, it is enough to prove that

\begin{eqnarray}
\label{s1}
&Q_{\boldsymbol \eta}[(r_{\kappa_1+\cdot}-r_{\kappa_1},\tau_{-r_{\kappa_1}}\eta'({\kappa_1}+\cdot))
\in F, A,\kappa_1=\tau_n]\\
\nonumber
&=
Q_{\delta_{0}}[
(r_\cdot,\eta(\cdot))\in F|D_0=\infty]
Q_{\delta_{0}}[A,\kappa_1=\tau_n],
\end{eqnarray}
for each natural $n$ and $A\in{\mathcal F}_{\tau_n}$.
Now, note that $\{\kappa_1=\tau_n\}=
\{D_n=\infty\}\cap\{D_1<\infty\}\ldots\cap\{D_{n-1}<\infty\}$.
But note that
for every $1\le j<n$, $\{D_n=\infty\}\cap\{D_j<\infty\}=\{D_n=\infty\}\cap\{D_j<\tau_n\}$.
Therefore

$$
\{\kappa_1=\tau_n\}=
\{D_n=\infty\}\cap\{D_1<\tau_n\}\ldots\cap\{D_{n-1}<\tau_n\}.
$$
Defining $B:=\{D_1<\tau_n\}\ldots\cap\{D_{n-1}<\tau_n\}$,
it follows that the right-hand side of (\ref{s1}) equals

\begin{eqnarray*}
&Q_{\boldsymbol \eta}[(r_{\kappa_1+\cdot}-r_{\kappa_1},\tau_{-r_{\kappa_1}}\eta'({\kappa_1}+\cdot))
\in F, A,\kappa_1=\tau_n]\\
&=
Q_{\boldsymbol \eta}[(r_{\tau_n+\cdot}-r_{\tau_n},\tau_{-r_{\tau_n}}
\eta'({\tau_n}+\cdot))
\in F, A,B, D_n=\infty]\\
&\!\!\! =
Q_{\boldsymbol \eta}[A,B,
Q_{\boldsymbol\eta}[(r_{\tau_n+\cdot}-r_{\tau_n},\tau_{-r_{\tau_n}}
\eta'({\tau_n}+\cdot))
\in F, D_n=\infty|{\mathcal F}_{\tau_n}]].
\end{eqnarray*}
Now, by translation invariance and the strong Markov property,
$
Q_{\boldsymbol\eta}[(r_{\tau_n+\cdot}-r_{\tau_n},\tau_{-r_{\tau_n}}
\eta'({\tau_n}+\cdot))
\in F, D_n=\infty|{\mathcal F}_{\tau_n}]
=Q_{\delta_0}[(r_{\cdot},
\eta(\cdot))
\in F, D_0=\infty]
$. Hence,

\begin{eqnarray}
\nonumber
&Q_{\boldsymbol \eta}[(r_{\kappa_1+\cdot}-r_{\kappa_1},\tau_{-r_{\kappa_1}}\eta'({\kappa_1}+\cdot))
\in F, A,\kappa_1=\tau_n]\\
\label{end}
&=
Q_{\boldsymbol \eta}[A,D_1<\tau_n,\ldots,D_{n-1}<\tau_n]
Q_{\delta_0}[(r_{\cdot},
\eta(\cdot))
\in F, D_0=\infty].
\end{eqnarray}
Choosing $F=\Omega'$ in the above equality we see that,
$Q_{\boldsymbol \eta}[A,\kappa_1=\tau_n]=
Q_{\boldsymbol \eta}[A,D_1<\tau_n,\ldots,D_{n-1}<\tau_n]
Q_{\delta_0}[D_0=\infty]$. Substituting this back into (\ref{end})
we conclude the proof of the proposition.
\end{proof}

\end{document}